\documentclass{amsart}
 
\usepackage{amscd,amssymb}

\theoremstyle{plain}
\newtheorem{thm}{Theorem}[section]

\theoremstyle{remark}
\newtheorem{rem}[thm]{Remark}

\renewcommand{\H}{\mathbb{H}\mkern1mu}

\newcommand{\bmat}{\left[\begin{smallmatrix}}
\newcommand{\emat}{\end{smallmatrix}\right]}
\newcommand{\bsmat}{\bigl[\begin{smallmatrix}}
\newcommand{\esmat}{\end{smallmatrix}\bigr]}
\newcommand{\Bsmat}{\Bigl[\begin{smallmatrix}}
\newcommand{\Esmat}{\end{smallmatrix}\Bigr]}
\newcommand{\bbsmat}{\biggl[\begin{smallmatrix}}
\newcommand{\eesmat}{\end{smallmatrix}\biggr]}
\newcommand{\BBsmat}{\Biggl[\begin{smallmatrix}}
\newcommand{\EEsmat}{\end{smallmatrix}\Biggr]}
\begin{document}
\title[Equivariant homotopy and deformations of diffeomorphisms]
{Equivariant homotopy and deformations of diffeomorphisms}
\author{C. Dur\'an}
\address{IMECC-UNICAMP, c.p. 6065, 13083-970, Campinas, SP, BRAZIL }
\email{cduran@ime.unicamp.br}        
\author{A. Rigas}
\address{IMECC-UNICAMP, c.p. 6065, 13083-970, Campinas, SP, BRAZIL }
\email{rigas@ime.unicamp.br}        

\begin{abstract}

We present a way of constructing and deforming diffeomorphisms of manifolds endowed with a Lie group action. 
This is applied to the study of exotic diffeomorphisms and involutions of spheres and to the equivariant homotopy of Lie groups. 
\end{abstract}


\maketitle

%
%
\section{Introduction}

In this paper we investigate the geometric and algebraic relation between two areas:

On the one hand, 
the study of the symmetries of geometric generators of homotopy groups:
it has long been a theme in homotopy theory to produce generators and elements 
of homotopy groups that are ``nice" with respect to symmetry and geometric properties 
(e.g. \cite{bott}, and, for more recent results,  \cite{thomashab} and the references therein).
On the other hand, we have the 
construction of deformations of diffeomorphisms, or the non-existence of such, i.e., orientation-preserving 
diffeomorphisms in different isotopy clasess, which give rise for example to 
exotic spheres (\cite{kervaire-milnor,duran}). Let us note that both areas 
have important apllications in physics, 
\cite{KR} for geometric generators of homotopy groups, and \cite{AB} for exotic 
phenomena.

We noticed the link between the two in a particular example, the relationship between a 
distiguished generator of $\pi_6(S^3)$ and an exotic diffeomorphism of $S^6$. In this 
paper we abstract this principle and reach what we call 
  ``equivariant J-process'', since it is analogous to the construction of the J-homomorphism 
in topology; we want to remark right away that in the present work we are mainly interested 
in the {\em algebra} of the J-process, instead of  the topology.   Our main result says that 
{\em equivariantly} symmetry-preserving deformations of elements of homotopy groups provide in a 
canonical way 
deformations of diffeomorphisms; and we give examples in which the lack of equivariance 
of the homotopy deformations implies that the deformed maps cease to be diffeomorphisms at some 
point (or, conversely, that no equivariant deformation exists). This process begins 
the distillation and abstraction of 
the phenomena that appears in the authors' research in exotic maps and involutions 
(\cite{duran,duran-mendoza-rigas, involutions, dupuri-matcont}), in order to search for a constructive and algebraic 
theory of exotic phenomena. 

The paper is organized as follows: in section \ref{main} we describe the main technique, 
which  is actually quite easy to prove; we believe that the relevant issue here is the 
abstraction of the principle and the consequences of its application in concrete  
examples. In section \ref{ourdiffeo} we explain 
the example that in fact motivated 
the main 
result: the equivariant differential geometry of a {\em Blakers-Massey} element (a generator of $\pi_6(S^3)$) and its relation to the 
exotic diffeomorphisms (i.e., not deformable to the identity through diffeomorphisms) of the 6-sphere, and exotic involutions of the 6-sphere and 
5-sphere. This setting will provide several applications of the main results, 
first deforming this exotic diffeomorphism of $S^6$ to a 
rational (still exotic) diffeomorphism, deforming the exotic involutions to rational maps, and then showing the non-existence of equivariant deformations of maps 
even though it is known that non-equivariant deformations exist. These applications are done in section \ref{applications}. 
We finally comment on some future directions in section \ref{CR}.

\medskip

\noindent{\bf Acknowledgments.} The authors would like to thank T. P\"uttmann for helpful discussions.

\section{The equivariant J- 
process} \label{main}

Let us describe the  equivariant J-process, which has two ingredients:

\begin{itemize}
 \item A Lie group $G$  acting differentiably on a manifold $M$ (from now on we assume that all groups, actions, maps are in the differentiable category). We denote the action of $G$ on $M$ by a dot $g\cdot m$, $g \in G, m\in M$.

\item A map  
$\alpha: M \to G$.
Note that in the case where $M$ is a $n$-sphere, the homotopy class of the map $\alpha$ then represents an element in the homotopy group $\pi_n(G)$.

\end{itemize}

   In this situation, we have

\medskip 

\noindent{\bf Definition.}
Define $J_\alpha: M \to M$, the J-process self map of $M$ associated to $\alpha$,  by  
$J_\alpha(m) = \alpha(m)\cdot m$. 

\medskip

We have then 

\begin{thm}\label{mainthm}
Let $J_\alpha: M \to M$ be a J-process associated to a map $\alpha: M \to G$ and an action of $G$ on $M$.
If $\alpha$ is $G$-equivariant with respect to the conjugation action  of $G$ on itself, then $J_\alpha$ is a bijection with   inverse 
 $(J_\alpha)^{-1} =  J_{\alpha^{-1}}$ and the powers are given by
${J_\alpha}^n = J_\alpha \circ \dots \circ J_\alpha = J_{\alpha^n}$

\end{thm}

\medskip

\begin{proof}

The equivariance hypothesis on $\alpha$ translates to $\alpha(g\cdot m) = g\alpha(m)g^{-1}$. Then, we just compute: 
\begin{eqnarray*}
 J_{\alpha^{-1}}(J_\alpha(m)) &=& J_{\alpha^{-1}}(\alpha(m) \cdot m) = 
\alpha^{-1}(\alpha(m)\cdot m)\cdot (\alpha(m)\cdot m)  \\
&=& ((\alpha(m)\alpha^{-1}(m)\alpha^{-1}(m))\cdot (\alpha(m) \cdot m) \quad  \text{by equivariance} \\
&=& \alpha^{-1}(m) \cdot (\alpha(m) \cdot m) = m \, . 
\end{eqnarray*}
The proof of the power formula is similar. 

\end{proof}

Consider now a one parameter family $\alpha_t: M \to G$ of differentiable maps. If this deformation satisfies 
the equivariance property for all $t$, then there is a one parameter family $J_{\alpha_t}$ of J-processes of $M$; theorem \ref{mainthm}  then guarantees that this deformation is {\em through diffeomorphisms}. We shall see 
in section \ref{applications} that sometimes a deformation through J-processes is guarateed to fail, that is,
 the map $J_{\alpha_t}$ must cease to be a diffeomorphism at some point in the parameter $t$ (of course in this 
 case the deformation is not equivariant).

We now endow the J-process construction with an additional structure: suppose  that in addition to the data in theorem
\ref{mainthm}, we have an involution $\delta: M \to M$. Then we have

\begin{thm} \label{mainthm-inv}
Let $M$ be a $G$-manifold, $\alpha: M\to G$ satisfying the hypothesis of theorem \ref{mainthm}. If in addition 
there is an involution $\delta$ of $M$ such that $\alpha(\delta(m)) = \alpha^{-1}(m)$, and $\delta$ commutes with the 
$G$-action (thus producing a $G\times \mathbb Z_2$-action on $M$),
then the J-process $\delta\circ J_\alpha$ is another involution of $M$. 
 \end{thm}
 
 \begin{proof}
 Compute: 
 \begin{eqnarray*}
 \delta J_\alpha (\delta J_\alpha (m)) &=& \delta J_\alpha (\delta (\alpha(m)\cdot m)) \\
 &=& \delta [\alpha(\delta (\alpha(m) \cdot m)) \cdot (\delta \alpha(m) \cdot m)]\\
 &=& \delta [\alpha^{-1}(\alpha(m) \cdot m) \cdot (\delta \alpha(m) \cdot m)], \text{ by the equivariance of $\delta$ and the inverse map,} \\
 &=& \delta [\alpha^{-1}(m) \cdot (\delta \alpha(m) \cdot m)], 
\text{ by the $G$-equivariance of $\alpha$,} \\
 &=& \delta \delta[\alpha^{-1}(m) \cdot ( \alpha(m) \cdot m)], 
\text{ since $\alpha$ and $\delta$ commute ,} \\
 &=& m, \text{ by $\delta^2=1$ and the group action property.}
\end{eqnarray*}
 \end{proof}

In the next sections we apply these results to an important special case.

\section{Blakers-Massey elements and exotic diffeomorphisms} \label{ourdiffeo}

\subsection{An equivariant generator of $\pi_6(S^3)$}

The objective of this section is to study   the geometry --in particular, the equivariant 
geometry-- of a distinguished generator of $\pi_6(S^3) \cong \mathbb Z_{12}$.  A map $f: S^6 \to S^3$   
is called by us a {\em Blakers-Massey element} if its class $[f] \in \pi_6(S^3)$ generates 
$\pi_6(S^3)$.  Recall that $\pi_6(S^3)$ 
classifies the principal $S^3$-bundles
over $S^7$, each element of $\pi_6(S^3)$ corresponding to an equivalence class of 
$S^3$-bundles over $S^7$. For example, the bundle 
$S^3 \cdots Sp(2) \to S^7$ corresponds to a generator of $\pi_6(S^3)$; the Gromoll-Meyer exotic sphere 
(\cite{gromoll-meyer})
of non-negative curvature is a given as a Riemannian quotient of the total space.  Also, Grove and Ziller 
(\cite{grove-ziller})
constructed cohomogeneity one metrics of non-negative curvature on principal $SO(4)$-bundles over $S^4$,   allowing the construction of metrics of non-negative curvature on the exotic 7-spheres that are bundles over $S^4$; some of these bundles are covered by
$S^3$-bundles over 
$S^7$.

It is a classic result \cite{james,toda} that $\pi_6(S^3) \cong \mathbb Z_{12}$; the traditional generator 
is obtained from the commutator of quaternions as follows: consider the map $\hat h:S^3 \times S^3 \to S^3$ given 
by $h(x,y) = xyx^{-1}y^{-1}$, where we consider $S^3$ as the group of unit quaternions. Since $h(1,y) = h(x,1) =1$, 
$h$ factors through a map $h:S^3 \wedge S^3 \cong S^6 \to S^3$ which generates $\pi_6(S^3)$ (\cite{toda}). In 
\cite{duran-mendoza-rigas}, by means of studying the geometry of geodesics of certain metrics on bundles over 
$S^7$, there
is a construction of a differentiable generator $b$ of $\pi^6(S^3)$ that is directly defined on $S^6$: 
consider 
the map $b: S^6 \to S^3$ defined as follows: the sphere $S^6$ is expressed as the set 
\[
S^6 = \{ (p,w) \in \H \times \H \,\, / \,\, \Re(p) = 0, |p|^2 + |w|^2 = 1 \} \, ,
\]
where $\H$ denotes the quaternions. Define $b:S^6 \to S^3$  by 
\[
b(p,w) = 
\begin{cases}
\frac{w}{|w|} e^{\pi p} \frac{\bar w}{|w|}, \quad w \neq 0 \\
-1, \quad \quad w = 0
\end{cases}
\]
where $e^x = \cos(|x|) + \sin(|x|)\frac{x}{|x|}$ denotes the exponential map of the Lie group $S^3$ of 
unit quaternions.  
The map $b$, which is a priori not even continuous at points of the form $(p,0)$, is in fact analytic, 
and it generates $\pi_6(S^3)$ 
(see \cite{duran-mendoza-rigas} and section \ref{deform}).

Having a differentiable map as a representative, obtained by geometrical methods, invites the study of its geometry, 
in particular with respect to group actions; we shall presently see that this distinguished generator has many symmetry properties and is a cohomogeneity one map.  
Write $SO(4) = S^3 \times S^3 /  \cong$, where $(q,r) \cong (-q,-r)$, and $SO(3) = S^3 / \pm 1$. 
 
 We represent $SO(4)$ and $SO(3)$ as follows 
Define $(q,r)\cdot (p,w) = (q p \bar q, r w \bar q)$. This representation embeds $SO(4)$ 
in the exceptional Lie group $G_2$; the group $SO(3)$ is represented  by the standard quaternionic conjugation 
$C_q(x) = qx\bar q$. The map $(q,r) \mapsto \pm q$ provides an epimorphism 
$\phi:SO(4) \to SO(3)$. Then, a simple computation shows that 
$b((q,r)\cdot (p,w))=  q b(p,w) \bar q$. Thus, we have the following commutative diagram, 
\[
\begin{array}{ccc}
SO(4) \times S^6 & \to & S^6 \\
(\phi,b) \downarrow &    & \downarrow b \\
SO(3) \times S^3 & \to & S^3 
\end{array}
\]

And thus the Blakers-Massey element $b$ is equivariant. Note that both actions are cohomogeneity one: 
the conjugation action on $S^3$ has as quotient the interval $[0,1]$, realized by the map $R:S^3 \to [-1,1]$, 
$R(\theta) = Re(\theta)$; the regular orbits are diffeomorphic to $S^2$ and $\pm 1$ are the two singular orbits.  
The $SO(4)$ action on $S^6$ has as the invariant function $S: S^6 \to \mathbb  [-1,1]$, 
$S(p,w) = |p|^2 - |w|^2$ ($= 2|p|^2 -1 = 1 - 2|w|^2$); the singular orbits are $S^2$ ($w=0$) and $S^3$ ($p=0$). The regular orbits are 
diffeomorphic to $S^2 \times S^3$, that is, the set $\{ (p,w) \,\, / \,\, |p| = r_1, |w|=r_2 \}$ when $r_1$ and $r_2$ 
are both non-zero. 

The following remark will be important in the sequel:

\begin{rem} \label{also_SO(3)_equivariant}
Consider the composition $\beta = P \circ b:S^6 \to SO(3)$, where $P$ is the canonical double cover projection $S^3 \to SO(3)$. 
Then taking the $\pm$ equivalence classes in the formula  $b((q,r)\cdot (p,w))=  q b(p,w) \bar q$ we also get a commutative 
diagram 
\[
\begin{array}{ccc}
SO(4) \times S^6 & \to & S^6 \\
(\phi,\beta) \downarrow &    & \downarrow \beta \\
SO(3) \times SO(3) & \to & SO(3) 
\end{array}
\]

\end{rem}
 
Let us also study at the powers of $b$:
\[
b^n(p,w) = 
\begin{cases}
\frac{w}{|w|} e^{n\pi p} \frac{\bar w}{|w|}, \quad w \neq 0 \\
(-1)^n, \quad \quad w = 0
\end{cases}
\]

Then 
all the powers of the Blakers-Massey element have the same equivariance properties. 
 
\subsection{The Blakers-Massey element and exotic maps}

This concrete map $b$ is a fundamental building block of exotic maps (degree one diffeomorphisms of spheres not isotopic to 
the identity, and free involutions not conjugate to the antipodal map);
see \cite{duran, duran-mendoza-rigas, involutions}.  Define 
\[
\sigma(p,w) = (b(p,w)\, p b(p,w)^{-1}, b(p,w)\, p b(p,w)^{-1})
\] 
then $\sigma:S^6 \to S^6$ is 
a degree one diffeomorphism not isotopic to the identity (\cite{duran,duran-mendoza-rigas}), and is 
a generator of the groups $\Gamma_7$ of isotopy classes of diffeomorphisms of $S^6$; (\cite{duran-mendoza-rigas}); 
it is knon (\cite{kervaire-milnor}) that 
$\Gamma_7 \cong \mathbb Z_{28}$.  
The map 
$\sigma^k$ represents $k \in \mathbb Z_{28}$, and thus in particular any diffeomorphism of $S^6$ can be deformed 
to  $\sigma^k$ for infinitely many $k$. Two maps $\sigma^k$ and $\sigma^\ell$ are isotopic to each other if 
and only if $k \equiv \ell \mod 28$, however let us remark that no explicit isotopy is known between 
$\sigma^k$ and $\sigma^{k + 28 r}$, $r \neq 0$. 
Note that the map $\sigma$ is $SO(3)$ 
(and not $SO(4)$) equivariant. 

It is not obvious that $\sigma$ is a diffeomorphism; in \cite{duran} 
this fact is established by indirect geometric methods, and in   \cite{duran-mendoza-rigas} it is 
observed that the inverse of $\sigma$ is given by 
$\sigma(p,w) = (b(p,w)^{-1}\, p b(p,w), b(p,w)^{-1}\, p b(p,w)$. Note that if $b$ were constant, 
this inverse is immediate; but since $b$ depends on $(p,w)$, in general such a result would not be true. 
However, we shall see that the equivariance 
is the structural reason for $\sigma$ being a diffeomorphism by applying the main theorem 
(compare 3.2 of \cite{duran-mendoza-rigas}); in order to do this, 
let us make precise how the J-process works in this case: 
let $B:S^6 \to SO(7)$ be given by $B(x) = \Delta \circ P \circ b(x) = \Delta \circ \beta$, where
$\Delta:SO(3) \to SO(7)$ is the diagonal embedding of $SO(3) \times SO(3)$ in $SO(7)$ (with a 1 in the middle),
i.e. 
\[
\Delta(T) = 
\begin{pmatrix}
T & 0 & 0 \\
0 & 1 & 0 \\
0 & 0 & T
\end{pmatrix}
\]
 
 Note that the image of $B$ falls into a subgroup of $SO(7)$ that is isomorphic to $SO(3)$. 
Then we are in the setup of theorem \ref{mainthm}: $B:S^6 \to SO(3)$ is equivariant with respect to the conjugation
action of $SO(3)$ in itself, and $\sigma$ can be simply written as the J-process 
$\sigma(x) = J_{B(x)}(x) = B(x) x$, since the projection of $S^3$ to $SO(3)$ is realized by the standard quaternionic conjugation 
$q \mapsto T_q$, $T_q(x) = qx \bar q$.  Thus we recover the results of \cite{duran-mendoza-rigas} in a more structural 
way: $\sigma$ is a diffeomorphism, and the powers of $\sigma$ are given by $\sigma^k(x) = B^k(x) x$. 

Also, $\sigma$ allows the construction of 
 {\em exotic involutions:} $(p,w) \mapsto -\sigma(p,w)$ is a free involution of $S^6$ that is not 
conjugate to the antipodal involution; since $\sigma$ acts on $(p,w)$ by a quaternionic conjugation, it preserves the real part of $w$ and thus this involution restricts to an (also exotic)
involution of the  $S^5$ defined by ${\rm Re}(w) = 0$, where 
it has a simple pictorial description (\cite{involutions,movie}). Again, the fact that $-\sigma$ is an involution 
follows immediately now by theorem \ref{mainthm-inv}, taking $\delta$ to be the antipodal involution of $S^n$ and 
$\alpha$ the Blakers-Massey element, since $b(-p,-w) = \overline{b(p,w)}$.

\section{Applications} \label{applications}

\subsection{Deformations of diffeomorphisms and involutions} \label{deform}

As mentioned in section \ref{main}, the equivariant J-process technique provides 
a ``canonical" way of deforming   diffeomorphisms $J_\alpha(m)$ through 
{\em equivariant} deformations of the corresponding maps $\alpha$. We will use this technique to deform 
the diffeomorphisms $\sigma^k: S^6 \to S^6$, which represent any isotopy class of diffeomorphisms of 
$S^n$, to  {\em rational} 
maps; this could pave the way to the algebro-geometric study of such exotic diffeomorphisms.   

The steps in the deformation are as follows: we will construct a homotopy $H(s,p,w):[0,1]\times S^6 \to S^6$ between 
the Blakers-Massey element and a rational map; this homotopy will be equivariant for all values of the deformation parameter 
$s$. 

By theorem \ref{mainthm}, the maps $J_{H_s}:S^6 \to S^6$ will be diffeomorphisms; and therefore this procedure 
furnishes a deformation between the exotic diffeomorphism $\sigma$ and a rational map $R$, in the same isotopy class 
of $\sigma$ and therefore also a generator of the group $\Gamma_7$. Then powers of $\sigma$ representing all other isotopy 
classes are also taken care of by 
theorem \ref{mainthm}, since $J_{H_s^k} = J_{H_s}^k = R^k$. 

In order to construct this deformation $H_s$, all we need to do is to reconsider proposition 1 of 
\cite{duran-mendoza-rigas} carefully and equivariantly.  Spelling out the exponential in the Blakers-Massey element, 
we have 
\[
b(p,w) = \cos(\pi |p|) + \frac{\sin (\pi |p|)}{|p|(1-|p|^2)} w p \bar w  \, .
\]
The functions $x\mapsto \sin(\pi x)$ and $x\mapsto x(1-x^2)$ are both odd, positive on $(0,1)$ and  have a zero of order 1 at 
$x=0$ and $x=1$; therefore $g(x) = \sin (\pi x)/(x(1-x^2))$ is an even, positive, differentiable function on $[0,1]$, 
(in particular this explains the analiticity of $b$).  We will homotop $g$ affinely to the constant function 1; in order 
to deal with $\cos(\pi |p|)$, we use the function $c(x) = 1 - 4x^2$. The function $c$ is the simplest even 
function satisfying the property that it has the same sign as $\cos(\pi x)$ on $[0,1]$. 

Now consider the $r(p,w) = (1-4|p|^2 + w p \bar w)$ and the  affine homotopy 
\[
\hat H(s,p,w) = \hat H_s(p,w) = (1-s) b(p,w) + s r(p,w) \, .
\]

For any $s$, the map $\hat H_s$ is equivariant with respect to conjugation since $b$ is equivariant and $r$ is a polynomial 
in $|p|,p$ and $w$. Also, the expressions in $|p|$ are all even, and therefore these maps can be written in terms 
of $|p|^2, p , w$, but since $p$ is purely imaginary, $|p|^2 = -p^2$ and all expressions involved are analytic 
expressions in the non-commuting quaternionic variables $p$, $w$ and $\bar w$.  

Rewriting $\hat H$, we have 
\[
\hat H_s(p,w) = [(1-s)\cos(\pi|p|) + s c(|p|)] + [(1-s) g(|p|) + s]w p \bar w \, .
\]

Then the sign properties of $c(x)$ and $g(x)$ imply that $\hat H(s,p,w)$ is never zero. Then 
$H(s,p,w) = \hat H(s,p,w)/|\hat H (s,p,w)|$ furnishes an equivariant homotopy between the Blakers-Massey element 
and  the map $Q: S^6 \to S^3$.
\[
Q(p,w) = \frac{1+4p^2 + wp\bar w}{\sqrt{(1+4p^2)^2 - |w|^4p^2}}
\]
 
 Now the map $R(p,w) = (Q(p,w) p \bar Q(p,w), Q(p,w) p \bar Q(p,w))$ is a rational diffeomorphism of $S^6$ that is not isotopic to the identity; its powers 
 are rational diffeomorphisms representing all isotopy classes of diffeomorphisms of $S^6$. Writing $R$ explicitly, we 
 have 
 \[
 R(p,w)=  
\left(\frac{(1+4p^2 + wp\bar w) p (1+4p^2 - wp\bar w)}{(1+4p^2)^2 - |w|^4p^2 }, 
\frac{(1+4p^2 + wp\bar w) w (1+4p^2 - wp\bar w)}{(1+4p^2)^2 - |w|^4p^2 } \right) \, .
\]

For each value $s$ of the deformation parameter, the map $H_s$ satisfies the hypothesis of theorem \ref{mainthm-inv} 
with respect to the antipodal involution of $S^n$. Therefore $-J_{H_s}$ is a deformation of the exotic 
involution $\sigma$ of $S^6$, through involutions that  also restrict to $S^5$, and, since close enough involutions are easily seen to be conjugate, these involutionas are all exotic.  At the end of the deformation 
we reach the involution $-R(p,w)$, which is a rational involution of $S^6$. Note that when restricted to $S^5$, 
defined by ${\rm Re}(w) = 0$,  $\bar w = -w$ and the map 
\[
 R(p,w)=  
-\left(\frac{(1+4p^2 - wpw) p (1+4p^2 + wpw)}{(1+4p^2)^2 - w^4p^2 }, 
\frac{(1+4p^2 - wp w) w (1+4p^2 + wpw)}{(1+4p^2)^2 - w^4p^2 } \right) \, .
\]
is an exotic involution of $S^5$ defined by a rational map in the non-commuting quaternionic variables $p$ and $w$.

Let us remark that deforming $\alpha$ through 
plain (not necessarily equivariant) homotopies produces a deformation of $J_\alpha$ that is not necessarily through 
diffeomorphisms.  We shall take advantage of this in the next section.

\subsection{The equivariant Serre problem} \label{sec-no-equi-serre}
It is know \cite{james,toda} that $\pi_6(S^3) \cong \mathbb Z_{12}$. However, no explict deformation of twelve times a generator 
is known. 
The authors call this the {\em Serre problem}: to find an 
{\em explicit} homotopy between the 12th power of the Blakers-Massey element and the identity, or, in other terms, to understand
how the quaternions are homotopy commutative in the 12 power.  
This problem is still open,
altough significant advances have been made recently (\cite{thomashab}); a solution to the Serre problem has far-reaching consequences, for example, the writing of explicit 
non-cancellation phenomena and new models for exotic spheres (see, e.g., \cite{rigas}). 
We show that, altough the generator and all its powers are represented by equivariant maps, 
{\em there is no equivariant solution} to the Serre problem (cf Theorem \ref{not-equi-serre}). We believe that, in addition to the statement of the theorem (which can probably be proven 
using standard methods of equivariant homotopy theory), the method of proof by using the relationship between 
equivariant homotopy and isotopy through explicit formulas  is of independent interest:
an equivariant homotopy would imply that the order of the group of homotopy 7-spheres divides 12 and we know 
it is isomorphic to $\mathbb Z_{28}$ (\cite{kervaire-milnor,eels-kuiper}); therefore such a homotopy is not possible. 
 We also extend this result to 
homotopies of the Blakers-Massey element inside other groups.
  
\medskip

\begin{thm} \label{not-equi-serre}
 There exists no differentiable homotopy $\phi:[0,1] \times S^6 \to S^3$ between 
$b^{12}$ and a constant map such that for each $t\in [0,1]$, $\phi(t,\cdot): S^6 \to S^3$ is $SO(3)$-equivariant.
\end{thm}

\begin{proof} 
We first adapt the Blakers-Massey element to the proposition above, by considering 
lifting to $S^3$ and considering $b: S^6 \to S^3$ as an $S^3$-equivariant map. An $SO(3)$ equivariant 
homotopy between $b^{12}$ and the constant map then lifts to an $S^3$-equivariant homotopy $b_t(p,w)$ 
such that $b_0(p,w) = 1$ and $b_1(p,w) = b(p,w)$. If such a homotopy exits, the maps 
\[
\sigma_t^{12}(p,w) = (b_t(p,w)^{12}\, p b_t(p,w)^{-12}, b_t(p,w)^{12}\, p b_t(p,w)^{-12})\, ,
\]
furnish an isotopy 
between $\sigma^{12}$ and the identity diffeomorphism. But by standard differential topology methods 
(see, for example, \cite{kosinski}), 
the map $\phi: \pi_0 (\mathrm{Diff}(S^6)) \to \Gamma_7$, $\phi(f) = D^7 \cup_f D^7$, is an isomorphism (\cite{milnor}), where 
$\Gamma_7$ is the group of differentiable structures on $S^7$ under the connected sum operation. Thus 
$\sigma^{12}$ represents 12 in the group  $\Gamma_7 \equiv  \mathbb Z_{28}$ (\cite{kervaire-milnor}), and we get a 
contradiction.  
\end{proof}

Note that the image of $B^k$ is contained in the chain of inclusions 
$SO(3) \subset SU(3) \subset G_2 \subset SO(7)$. Theorem \ref{not-equi-serre} states that, even though 
$B^{12}$ is homotopic (inside of $SO(3)$) to the constant map, no equivariant homotopy can exist. Now the 
``Serre problem" for all the other groups has been solved: there exist explicit generators $\gamma$ of $
\pi_6(SU(3)) \cong \mathbb Z_ 6$, $\delta$ of $\pi_6(G_2) \cong \mathbb Z_3$ and explicit homotopies 
between $\gamma^6$ and the constant map \cite{puri1,thomashab} and $\delta^3$ and the constant map 
\cite{rigas0}.
Also $\pi_6(SO(7)) = 0$. 
Explicit homotopies between $b$ and $\gamma$, $\delta$ and the constant map inside of the respective groups can be  
constructed using the geometry of the chain
$SO(3) \subset SU(3) \subset G_2 \subset SO(7)$ \cite{thomashab}. Thus, not only  
we have that $b^6, b^3$ and $b$ are homotopic to the constant map inside of $SU(3), G_2$ and $SO(7)$; 
explicit homotopies can be written.

\begin{thm} \label{no-equi-other-serre}
The  maps $b^6, b^3, b$ are homotopic to the constant map  in  $SU(3), G_2, SO(7)$,  
respectively, through explicit homotopies. However, no $SO(3)$-equivariant homotopy exists.
\end{thm}

\begin{proof} 
All the groups in the chain $SO(3) \subset SU(3) \subset G_2 \subset SO(7)$ act on 
$S^6$ through the canonical action of $SO(7)$. Mimic the proof of Theorem 1 with $\sigma^6, \sigma^3$ and 
$\sigma$ in place of $\sigma^{12}$.   
\end{proof} 

The symmetry-breaking mechanism of these homotopies is beautifully illustrated in the structure of the group $\pi_6(SU(3))$
\cite{puri1}: one way of determining the structure 
of the homotopy of topological groups is by finding a generator $A$ such that $A^k= e$, the identity of the group; this 
is the way that this was done for $\pi_6(G_2) \cong \mathbb Z_3$, by finding a generator such that $A^3 = e$ (\cite{rigas0}). 
The power map $A \mapsto A^k$ of matrices is of course equivariant under conjugation and this process provides equivariant 
deformations.  However, in the case of $SU(3)$, the 
homotopy uses a {\em deformation} of the product of matrices through the Cartan subalgebra of $SU(3)$ and the 
symmetry is broken; see \cite{puri1} for details.

\section{Concluding remarks} \label{CR}

First we want to note that all the constructions in \cite{duran-mendoza-rigas,involutions} 
related to the Blakers-Massey element and the exotic diffeomorphisms 
can be generalized by substituting all the quaternions involved by Cayley numbers and modifying the relevant dimensions. 
However, this passage involves non-trivial modifications of the techniques in the proofs. The unit quaternions are a group, and thus 
the equivariance properties make sense; the unit Cayley numbers are not a group and what is the right extension 
of theorem \ref{mainthm} and its applications remains to be seen.

The construction in the main theorem, the applications given and the computations in section 3 of \cite{duran-mendoza-rigas} 
suggests that there exists an ``algebra of exoticity", which is yet to be described.

It would also be interesting to follow the known, non-equivariant homotopies of the respective powers of the 
Blakers-Massey element  and the identity in $SU(3)$, $G_2$ and $SO(7)$  
and study the   associated J-process self-maps of $M$, which at some point must cease to be diffeomorphisms 
and determine the structure of the singularities that appear, which could shed some light in the general 
question of what makes a diffeomorphism exotic, or, better, how does one detect the exoticity of a degree 
one diffeomorphism that is given by a formula. 

%

\end{document}